\documentclass[12pt,a4paper,reqno]{amsart}
\usepackage{amsfonts}
\usepackage{amsthm}
\usepackage{amsmath}
\usepackage{amssymb}
\usepackage{amscd}
\usepackage[latin2]{inputenc}
\usepackage{t1enc}
\usepackage[mathscr]{eucal}
\usepackage{indentfirst}
\usepackage{graphicx}
\usepackage{graphics}
\usepackage{pict2e}
\usepackage{epic}
\numberwithin{equation}{section}
\usepackage[margin=2.9cm]{geometry}
\usepackage{epstopdf}
\usepackage{fancyhdr}
\usepackage[colorlinks]{hyperref}

\newtheoremstyle{break}
  {\topsep}{\topsep}%
  {\normalfont}{}%
  {\bfseries}{.}%
  {\newline}{}%
\theoremstyle{plain}
\newtheorem{teo}{Theorem}[]
\newtheorem*{teointro}{Theorem}
\theoremstyle{break}
\newtheorem{teo-bis}[teo]{Theorem}
\theoremstyle{break}
\newtheorem{defin-bis}[teo]{Definition}
\theoremstyle{plain}       %definizione ambiente teorema a capo
\newtheorem{prop}[teo]{Proposition}    %definizione ambiente proposizione
\newtheorem{lem}[teo]{Lemma}            %definizione ambiente lemma
\theoremstyle{definition}               %stile roman
\newtheorem{defin}[teo]{Definition}%definizione ambiente definizione
\theoremstyle{definition}                   %stile per osservazioni

\theoremstyle{remark}

\theoremstyle{remark}
\newtheorem{remk}[teo]{Remark} 
\theoremstyle{remark}
\newtheorem{remk-count}[teo]{Remark}
\makeatletter
\newcommand\thankssymb[1]{\textsuperscript{\@fnsymbol{#1}}}
\makeatother

\begin{document}

\title[]{Representations of automorphism groups\\on the homology of matroids}

\author[]{Luca Moci\thankssymb{1}\thankssymb{2}}

\author[]{Gian Marco Pezzoli\thankssymb{1}}
\thanks{\thankssymb{1} Universit\`a di Bologna}
\thanks{\thankssymb{2} supported by PRIN 2017YRA3LK}
       
        %Gian Marco Pezzoli%\thanks{
                     % Universit\`a di Bologna
                  %}
        %}

%%%%%%%%%%%%%%%%%%%%%%%%%%%
% abstract, keywords and Subject classification are optional.
%%%%%%%%%%%%%%%%%%%%%%%%%%%
\begin{abstract}
    Given a group $G$ of automorphisms of a matroid $M$, we describe the representations of $G$ on the homology of the independence complex of the dual matroid $M^*$. These representations are related to the homology of the lattice of flats of $M$, and (when $M$ is realizable) to the top cohomology of a hyperplane arrangement. Finally, we analyze in detail the case of the complete graph, which has applications to algebraic geometry.
\end{abstract}

\maketitle

\tableofcontents

\section{Introduction}
In the last years, matroid theory received increasing attention from geometers because of its multiple connections with algebraic geometry, algebraic topology and representation theory. One of the equivalent ways of defining a matroid $M$ is by specifying the family of its independent sets, which is an abstract simplicial complex $IN(M)$. A group $G$ of automorphism of $M$ acts naturally on $IN(M)$; this gives rise to representations of $G$ on the homology of $IN(M)$, which are studied in this article. 

One motivation for our work comes from a paper by de Cataldo, Heinloth and Migliorini, that computes the supports of the perverse cohomology sheaves of the Hitchin fibration for $GL_m$ over the locus of reduced spectral curves, studying the related Cattani-Kaplan-Schmid complex (\cite{DHM}). 
The dual graph of such a spectral curve is the complete graph, and the action of the symmetric group on the irreducible components of the curve yields an action on its vertices,  hence on the independence complex of the dual matroid of the graph. A crucial step in the analysis performed in \cite{DHM} is then to determine the representations of the symmetric group on the homology of this independence complex. 

In this paper we describe representations on the homology of matroids in full generality. Namely, given a matroid $M$ of rank $r$ on $n$ elements and a group $G$ of automorphisms of $M$, we consider the independence complex $\hbox{IN}(M^{*})$ of the matroid dual to $M$. We prove that the reduced homology of this simplicial complex, up to a shift and to a sign, is isomorphic to the reduced homology of the non-spanning complex of $M$, $NS(M)$, and to the reduced homology of the order complex of the lattice of flats of $M$, $\mathcal L (M)$:
\begin{teointro}
The following representations of $G$ are isomorphic for every $i \geq 0$ (and nonzero only for $i=r-2$):
\begin{enumerate}
    \item\: $\widetilde{H}_{n-3-i}(\hbox{IN}\,(M^{*})) \otimes \textnormal{sgn}$
    \item\: $\widetilde{H}_{i}(\hbox{NS}\,(M))$
    \item\:$\widetilde{H}_{i}(\mathcal{L}\,(M))$
\end{enumerate}
where $n$ is the cardinality of the ground set of $M$ and \textnormal{sgn} is the sign representation (restricted from $\mathfrak{S}_n$ to its subgroup $G$).
\end{teointro}
Here, the isomorphism between $(1)$ and $(2)$ holds more generally for any simplicial complex, being a consequence of Alexander duality (see Theorem \ref{teo:fonda}). The proof that we give in Section \ref{section:alexander}  is inspired by that of Bj$\ddot{\hbox{o}}$rner and Tancer (\cite{Bj2}), but keeps track of the action of the group $G$. 

Also the isomorphism between $(2)$ and $(3)$ is a consequence of a more general phenomenon: indeed, in Section \ref{section:crosscut} we develop an equivariant version of Folkman's machinery of \emph{cross-cuts} \cite{folk}: see in particular our Theorem \ref{teo:crosscutfixed}.

In Section \ref{section:matroid} we specialize our results to the case of matroids, obtaining the above-mentioned isomorphisms.

Furthermore, if the matroid $M$ is realizable, then it is naturally associated with a hyperplane arrangement $\mathcal{A}$. The cohomology of the complement $\mathcal{C}(\mathcal{A})$ of the arrangement admits a well-known presentation in terms of $M$, due to Orlik and Solomon (\cite{orliksolomon}). In Section \ref{section:arrangement} we show that the top-degree part of this cohomology is isomorphic as a representation of $G$, up to a sign, to the reduced homology of the dual matroid $M^{*}_{\mathcal{A}}$ associated to $\mathcal A$:
$$H^{r}(\mathcal{C}(\mathcal{A}))  \simeq_{G} \widetilde{H}_{n-r-1} (IN(M^{*}_{\mathcal{A}})) \otimes \textnormal{sgn}$$
%where $r$ is the rank of the arrangement and $n$ is the cardinality of the ground set of $M_{\mathcal{A}}$ 
(see Theorem \ref{teo:iperpiani}). This statement is a consequence of our main theorem and of results of Orlik and Solomon \cite{orliksolomon}.

In Section \ref{section:kook} we focus on the case when $M$ is realized by a \textit{coned graph} or by a \textit{ complete bipartite graph}. Then our Theorem \ref{teo:combined}, combined with results by Kook and Lee (\cite{kook,kook2}), yields isomorphisms with the $G-$modules of \emph{edge-rooted} and \emph{B-edged rooted forests} (\ref{gra1}, \ref{gra2}).

Finally, in Section \ref{section:completegraph} we specialize our results to the case in which $M$ is the matroid of the complete graph $K_m$, or equivalently of the root system of type $A_{m-1}$, which is the case of interest in \cite{DHM}. For this matroid, whose group of automorphisms is the symmetric group  $\mathfrak{S}_m$ and whose lattice of flats is the partition lattice $\Pi_m$, a result of Stanley \cite{stan1} allows explicit determination of the representations:
$$
\widetilde{H}_{n-m} (IN(M^{*}(K_m))) \simeq_{\mathfrak{S}_m} \textnormal{ind}_{C_m}^{\mathfrak{S}_m}(e^{2\pi i / m})
$$
where $n= \binom{m}{2}$ is the number of edges $K_m$, and $C_m$ is the subgroup generated by an $m$-cycle in $\mathfrak{S}_m$ (see also Remark \ref{miglio}).

It is natural to wonder if a similar description can be provided for root systems of other types. This seems to be a hard question, though, being equivalent to the following conjecture of Lehrer and Solomon (\cite[Conj. 1.6]{lehrersolomon}):
$$
H^{p}(\mathcal{C}(\mathcal{A}_W)) \simeq_{W} \bigoplus_c \textnormal{Ind}^{W}_{Z(c)}(\xi_c) \qquad p=0,\dots,\textnormal{rank}(W).
$$

A different direction of research, that we hope to develop in future papers, is to explicitly describe the representation arising for different classes of graphs or matroids.
\vspace{0.2cm}
\paragraph {\bfseries Acknowledgements} We are grateful to Emanuele Delucchi, Luca Migliorini and Roberto Pagaria for inspiring discussions and very valuable suggestions. We also wish to thank the anonymous referees for many helpful remarks.

\section{Representations and Alexander duality}\label{section:alexander}
We recall here some basic facts in combinatorial topology. For more details the reader can refer to \cite{hatcher}.  Let $K$ be an abstract simplicial complex with vertex set $V$ with $|V|=n$.
For $\sigma \in K$, let 
$$
\overline{\sigma}=V \smallsetminus \sigma.
$$

\begin{defin}
The \textit{Alexander dual} of $K$ is the simplicial complex on the same vertex set defined by
$$
K^{*}=\{\sigma \subseteq V \ | \ \overline{\sigma} \notin K\}.
$$
It is easy to see that $K^{**}=K$.
\end{defin}
Let $G$ be a finite group of automorphisms of the face poset $(K,\subseteq)$. 
Then $G$ is a subgroup of the symmetric group $\mathfrak{S}_n$ on $V$, made out of the vertex maps
$
\begin{array}{cccc} 
g:& V & \longrightarrow & V \\
\end{array}
$
such that whenever the vertices $j_{1},\dots,j_{i+1}$ span an $i$-simplex of $K$, the points $g(j_{1}),\dots,g(j_{i+1})$ span an $i$-simplex of $K$. Therefore $g$ induces a simplicial homeomorphism $\tilde g$, and $\tilde g$ induces a chain-isomorphism
 $\tilde g_{\#}$ on the group of oriented $i$-chains in the following way:
 $$
\begin{array}{cccc} 
\tilde g_{\#,i}:& C_{i}(K,\mathbb{C}) & \longrightarrow & C_{i}(K,\mathbb{C}).\\
& [j_{1},\dots,j_{i+1}] & \longmapsto & [g(j_{1}),\dots,g(j_{i+1})]
\end{array}
 $$
 
Moreover $\tilde g_{\#}$ induces an isomorphism on the reduced homology groups $\widetilde{H}_{i}(K)$ (and the reduced cohomology groups $\widetilde{H}^{i}(K)$):
$$
\begin{array}{cccc} 
\rho_{i,g}:& \widetilde{H}_{i}(K,\mathbb{C}) & \longrightarrow & \widetilde{H}_{i}(K,\mathbb{C}).
\end{array} \qquad 
$$
This defines representations of $G$ on the $\mathbb{C}$-vector spaces $\widetilde{H}_{i}(K,\mathbb{C})$, i.e., homomorphisms
$$
\begin{array}{cccc} 
\rho_i:& G & \longrightarrow & \textnormal{GL}(\widetilde{H}_{i}(K,\mathbb{C})). \\
& g & \longmapsto & \rho_{i,g} 
\end{array}
$$
It follows from the definition of $K^{*}$ that $G$ is also a finite group of automorphisms of the \textit{face poset} of $K^{*}$. Therefore, following the construction above, we get representations $\rho^{*^i}$ of $G$ on the reduced cohomology of  $K^{*}$:
%forse asterisco piu indice in alto e' un po' pesante come notazione, sposterei l'indice in basso
$$
\begin{array}{cccc} 
\rho^{*^i}:& G & \longrightarrow & \textnormal{GL}(\widetilde{H}^{i}(K^{*},\mathbb{C})). \\
& g & \longmapsto & \rho^{* i}_{g} 
\end{array}
$$
\begin{teo}\label{teo:fonda}
	Let $K$ be an abstract simplicial complex and let $K^*$ be its Alexander dual. Let $G$ be a finite group of automorphisms of the \textit{face poset} of $K$. Then:
	$$
	\rho_i \simeq \rho^{*^{n-i-3}} \ \otimes \ \textnormal{sgn}
	$$
	where $n=|V|$ and \textnormal{sgn} is the sign representation (restricted from $\mathfrak{S}_n$ to its subgroup $G$). Or, equivalently, we have the following isomorphism of $\mathbb{C}[G]$-modules: 
	$$
	\widetilde{H}_{i}(K,\mathbb{C}) \simeq_{G} \widetilde{H}^{n-i-3}(K^{*},\mathbb{C})  \otimes  \textnormal{sgn}.
	$$
\end{teo}

Our proof follows from Bj$\ddot{\hbox{o}}$rner and Tancer (\cite{Bj2}), but carefully records the action $G$. We introduce some notations: 
let us denote by $\{1,2,\dots, n\}$ the elements of $V$.  For $j \in \sigma \in K$, we define the \textit{sign}
$$
\textnormal{sgn}(j,\sigma)=(-1)^{i-1}
$$ 
where $j$ is the $i$-th smallest element of the set $\sigma$, and 
$$
p(\sigma)=\prod_{j \in \sigma} (-1)^{j-1}.
$$
Take an $i$-face $\sigma=\{j_1,\dots,j_{i+1}\} \in K$ with $j_1 < \cdots < j_{i+1}$. We write $e_{\sigma}$ to denote the oriented simplex associated to $\sigma$ considered with an increasing order of its elements:
$$
e_{\sigma}=[j_1,\dots,j_{i+1}].
$$
For every $g \in G$, we denote 
$$
g.\sigma=\{g(j_{1}),\dots,g(j_{i+1})\} \quad \hbox{and} \quad g. e_{\sigma}=[g(j_{1}),\dots,g(j_{i+1})].
$$

The $g(j_{1}),\dots,g(j_{i+1})$ are not necessarily in ascending order: let $\tau \in \mathfrak{S}_{i+1} \subseteq \mathfrak{S}_n$ be the permutation that rearranges the elements in ascending order, and fixes the elements that are not in $g.\sigma$, so that $\tau.(g. e_{\sigma})=e_{g.\sigma}$. We also define:
$$
c(g,\sigma)=\textnormal{sgn}(\tau).
$$

Since $\tau^{-1}$ permutes the elements of $e_{g.\sigma}$ we obtain:
$$
g. e_{\sigma}=\tau^{-1}.(e_{g.\sigma})=\textnormal{sgn}(\tau^{-1})e_{g.\sigma}=\textnormal{sgn}(\tau)e_{g.\sigma}=c(g,\sigma) e_{g.\sigma}.
$$
Similarly, we define a permutation $\overline{\tau} \in \mathfrak{S}_{n-1-i} \subseteq \mathfrak{S}_n$ which rearranges the elements of $g. e_{\overline{\sigma}}$ in ascending order:
$$
\overline{\tau}.(g. e_{\overline{\sigma}})=e_{g.\overline{\sigma}} \qquad g. e_{\overline{\sigma}}=c(g,\overline{\sigma})e_{g.\overline{\sigma}}.
$$
We can now formulate an important lemma that will prove to be crucial for the proof of Lemma \ref{teo:trequattro}:
\begin{lem}\label{lem:vacci}
Let $V=\{1,\dots,n\}$ and let $\sigma	\subseteq V$. Then, for every $g \in \mathfrak{S}_n$, we have the following:
\begin{equation}\label{eq:res}
p(\sigma) \ \textnormal{sgn}(g) \ c(g, \overline{\sigma})=c(g,\sigma) \ p(g.\sigma).
\end{equation}
\end{lem}
\begin{proof}
For every $g \in \mathfrak{S}_n$, we define a permutation $g'=\overline{\tau} \tau g$. First we apply the permutation $g$ to $\sigma$ and $\overline{\sigma}$. Then applying $\tau$ and $\overline{\tau}$, we rearrange in ascending order both $g. e_{\sigma}$ and $g. e_{\overline{\sigma}}$. 

As we have defined it, $g'$ is a permutation of $\mathfrak{S}_n$ such that:
\begin{itemize}
\item[] if $i,j \in \sigma$ with $i<j$ then $g'(i)<g'(j)$ and
\item[] if $i,j \in \overline{\sigma}$ with $i<j$ then $g'(i)<g'(j)$.
\end{itemize}
In particular, we have that:
$$
g'. e_{\sigma}=e_{g'. \sigma} \qquad \hbox{and}\qquad g'. e_{\overline{\sigma}}=e_{g'. \overline{\sigma}}.
$$
We can express $g$ in the following way $g=\tau^{-1} \circ \overline{\tau}^{-1}\circ  g'$. It easy to see that $p(g.\sigma)=p(g'.\sigma)$. Thus, Equation (\ref{eq:res}) is equivalent to:
\begin{align}\label{eq:kik}
    \nonumber
    p(\sigma)  \textnormal{sgn}(\tau^{-1} \overline{\tau}^{-1} g')  c(g, \overline{\sigma})=c(g,\sigma)  p(g.\sigma) &\Leftrightarrow p(\sigma)  \textnormal{sgn}(\tau^{-1})\textnormal{sgn}(\overline{\tau}^{-1}) \textnormal{sgn}(g')  \textnormal{sgn}(\overline{\tau})=\textnormal{sgn}(\tau) p(g'.\sigma)\\ 
    \nonumber
    &\Leftrightarrow p(\sigma) \ \textnormal{sgn}(g')= p(g'.\sigma)\\
    \nonumber &\Leftrightarrow \Big( \prod_{i \in \sigma}  (-1)^{i-1} \Big)  \textnormal{sgn}(g')=\prod_{i \in \sigma}  (-1)^{g'(i)-1}\\
    \nonumber
    &\Leftrightarrow \Big( \prod_{i \in \sigma}  (-1)^{i-g'(i)} \Big) \textnormal{sgn}(g')=1\\
    &\Leftrightarrow \prod_{i \in \sigma}  (-1)^{i-g'(i)}=\textnormal{sgn}(g').
\end{align}
In order to prove Equation (\ref{eq:kik}), let $i \in \sigma$ be the $k$-th element of $e_{\sigma}$ and we define:
\begin{align*}
    A_i&=\{(i,j) \ | \ j \in \overline{\sigma}, \ i<j, \ g'(i)>g'(j)\} \hbox{ and}\\
    B_i&=\{(j,i) \ | \ j \in \overline{\sigma}, \ j<i, \ g'(j)>g'(i)\}.
\end{align*}

We have that:
$$
\textnormal{sgn}(g')=(-1)^{\sum_{i \in \sigma} (|A_i|+|B_i|)}.
$$
Let us assume that $i<g'(i)$. It is easy to see that $|B_i|=0$. Furthermore:
$$
|\{(i,j) \ | \ j \in \overline{\sigma},\ i<j\}|=(n-i)-(|\sigma|-k)=n-i-|\sigma|+k
$$
and
\begin{align*}
|\{(i,j) \ | \ j \in \overline{\sigma},\ i<j, \ g'(i)<g'(j)\}|&=|\{(i,j) \ | \ j \in \overline{\sigma},\ g'(i)<g'(j)\}|
\\
&=(n-g'(i))-(|\sigma|-k)=n-g'(i)-|\sigma|+k.
\end{align*}
By subtracting term by term the two equalities above, we get:
$$
|A_i|=(n-i-|\sigma|+k)-(n-g'(i)-|\sigma|+k)=g'(i)-i.
$$
Similarly, if $i>g'(i)$ we have that $|B_i|=i-g'(i)$ and $|A_i|=0$. Therefore:
$$
|A_i|+|B_i|=|g'(i)-i|.
$$
It follows that:
$$
\textnormal{sgn}(g')=(-1)^{\sum_{i \in \sigma} |g'(i)-i|}=\prod_{i \in \sigma}\ (-1)^{|g'(i)-i|}=\prod_{i \in \sigma} \ (-1)^{g'(i)-i}.
$$
\end{proof}

Let $2^{V}$ be the full simplex with vertex set $V$.
\begin{lem}\label{teo:trequattro}
Let $K$ be a simplicial complex with ground set $V$ of size $n$. Then 
$$
\widetilde{H}_{i+1}(2^{V},K) \simeq \widetilde{H}^{n-i-3}(K^{*}).
$$
Furthermore, if we consider the following representations of the group $G$ on the homology spaces of $(2^{V},K)$ and $K^*$:
$$
\begin{array}{cccc} 
\alpha_{i+1}:& G & \longrightarrow & \textnormal{GL}(\widetilde{H}_{i+1}(2^{V},K),\mathbb{C}))
\end{array}  \qquad \hbox{and}
$$
$$
\begin{array}{cccc} 
\rho^{*^{n-i-3}}:& G & \longrightarrow & \textnormal{GL}(\widetilde{H}^{n-i-3}(K^*,\mathbb{C}))
\end{array},
$$
we have that
$$
\alpha_{i+1} \simeq \rho^{*^{n-i-3}} \ \otimes \ \textnormal{sgn}.
$$
Or equivalently
$$
\widetilde{H}_{i+1}(2^{V},K) \simeq_{G} \widetilde{H}^{n-i-3}(K^{*}) \otimes \textnormal{sgn}.
$$
\end{lem}

\begin{proof}
The chain complex for reduced homology of the pair $(2^{V},K)$ is the complex:
$$
\cdots \mathcal{R}_{i+1} \xrightarrow{d_{i+1}} \ \mathcal{R}_i \xrightarrow{d_i} \ \mathcal{R}_{i-1} \xrightarrow{d_{i-1}} \ \cdots, \qquad i \in \mathbb{Z}
$$
where $\mathcal{R}_i=\langle e_{\sigma} \ | \ \sigma \subseteq V, \ \sigma \notin K, \ \textnormal{dim}(\sigma)=i \rangle$, and the $d_i$'s are the unique homomorphisms satisfying:
$$
d_{i}(e_{\sigma})= \underset{\sigma \smallsetminus k \notin K}{\sum_{k \in \sigma}} \ \textnormal{sgn}(k, \sigma) \ e_{\sigma \smallsetminus k}.
$$

The cochain complex for reduced cohomology of $K^{*}$ is the complex:
$$
\cdots \xrightarrow{\delta_{i-1}} \ C^{i-1} \xrightarrow{\delta_{i}} \ C^{i} \xrightarrow{\delta_{i+1}} \ \cdots, \qquad i \in \mathbb{Z}
$$
where
\begin{align*}
    C^{i}&=\langle e_{\sigma}^{*} \ | \ \sigma \subseteq V, \ \textnormal{dim}(\sigma)=i, \ \sigma \in K^{*} \rangle \\
    &=\langle e_{\sigma}^{*} \ | \ \sigma \subseteq V, \ \textnormal{dim}(\overline{\sigma})=n-i-2, \ \overline{\sigma} \notin K \rangle
\end{align*}
and the $\delta_i$'s are the unique homomorphisms satisfying:
$$
\delta_i(e^{*}_{\sigma})=\underset{\sigma \cup k \in K^{*}}{\sum_{k \notin \sigma}} \ \textnormal{sgn}(k, \sigma \cup k) e^{*}_{\sigma \cup k}=
\underset{\overline{\sigma} \smallsetminus k \notin K}{\sum_{k \in \overline{\sigma}}} \ \textnormal{sgn}(k, \sigma \cup k) e^{*}_{\overline{\overline{\sigma} \smallsetminus k}}.
$$

Let $\phi_i$ be the following isomorphism:
\begin{equation}\label{iso:phi}
\begin{array}{cccc} 
\phi_i:& \mathcal{R}_i & \longrightarrow & C^{n-i-2} \\ 
& e_{\sigma} & \longmapsto & p(\sigma) \ e^{*}_{\overline{\sigma}} 
\end{array} \quad \hbox{ for } \sigma \notin K \hbox{ with } \textnormal{dim}(\sigma)=i.
\end{equation}
We then have the following diagram:
$$
\begin{CD}
@>d_{i+1}>> \mathcal{R}_i @> d_i  >>  \mathcal{R}_{i-1}  @>d_{i-1}>> \\
@. @V\phi_{i}VV @V\phi_{i-1}VV\\
@>\delta_{n-i-2}>> C^{n-i-2} @ >\delta_{n-i-1} >> C^{n-i-1} @>\delta_{n-i}>>
\end{CD}
$$
\newline
We know from the proof of Lemma 4.2 of \cite{Bj2} that
\begin{equation}\label{eq:trova}
\phi_{i-1} \circ d_i=\delta_{n-i-1} \circ \phi_i.
\end{equation}
Thus, we have that
$$
\widetilde{H}_{i+1}(2^{V},K) \simeq \widetilde{H}^{n-i-3}(K^{*}).
$$
We now define the following two representations:

\begin{align*}
   &\begin{array}{cccc} 
    \rho_{1}:& G & \longrightarrow & \textnormal{GL}(\mathcal{R}_i) \\
    & g & \longmapsto & \rho^{1}_{g} 
    \end{array} \quad \hbox{by} \quad  \begin{array}{cccc} 
    \rho^{1}_g:& \mathcal{R}_i & \longrightarrow &\mathcal{R}_i \\
    & e_{\sigma} & \longmapsto & g.e_{\sigma}
    \end{array}
    \quad \hbox{and} \\
    &\begin{array}{cccc} 
    \rho_{2}:& G & \longrightarrow & \textnormal{GL}(C^{n-i-2}\otimes \mathbb{C}) \\
    & g & \longmapsto & \rho^{2}_{g} 
    \end{array} \quad \hbox{by} \quad \begin{array}{cccc} 
    \rho^{2}_g:& C^{n-i-2}\otimes \mathbb{C}& \longrightarrow &C^{n-i-2}\otimes \mathbb{C} \\
    & e^{*}_{\overline{\sigma}} \otimes 1 & \longmapsto & g.e^{*}_{\overline{\sigma}} \otimes \textnormal{sgn}(g)
    \end{array}
\end{align*}
for $\sigma \notin K$ with dim$(\sigma)=i$. We want to show that this two representations are isomorphic. We extend the isomorphism (\ref{iso:phi}):
$$
\begin{array}{cccc} 
\tilde{\phi}_i:& \mathcal{R}_i & \longrightarrow & C^{n-i-2}\otimes \mathbb{C} \\ 
& e_{\sigma} & \longmapsto & p(\sigma) \ e^{*}_{\overline{\sigma}} \otimes 1
\end{array} \quad \hbox{ for } \sigma \notin K \hbox{ with } \textnormal{dim}(\sigma)=i.
$$
To prove that $\rho_1 \simeq \rho_2$ we have to show that the following diagram commutes for every $g \in G$:
$$
\begin{CD}
\mathcal{R}_i @> \rho^1_{g}  >>  \mathcal{R}_{i-1}   \\
@VV\tilde{\phi}_{i}V @VV\tilde{\phi}_{i}V\\
 C^{n-i-2}\otimes \mathbb{C} @ >\rho^2_{g} >> C^{n-i-1} \otimes \mathbb{C}
\end{CD}
$$

We have to prove that the following equation holds:
\begin{equation}\label{eq:comma}
\rho^2_g \circ \tilde{\phi}_i=\tilde{\phi}_i \circ \rho^1_g.
\end{equation}

Note that for each $e_{\sigma} \in \mathcal{R}_i$,
\begin{align*}
    (\rho^2_g \circ \tilde{\phi}_i)(e_{\sigma})&=\rho^2_g(p(\sigma) \ e^{*}_{\overline{\sigma}} \otimes 1)=p(\sigma)  \ g.e^{*}_{\overline{\sigma}} \otimes \textnormal{sgn}(g)=
    p(\sigma)  \  \textnormal{sgn}(g) c(g,\overline{\sigma}) e^{*}_{g.\overline{\sigma}} \otimes 1, \hbox{ and}\\
    (\tilde{\phi}_i \circ \rho^1_g)(e_{\sigma})&=\tilde{\phi}_i(g.e_{\sigma})=
    \tilde{\phi}_i(c(g,\sigma)e_{g.\sigma})=p(g.\sigma)     c(g,\sigma)e^{*}_{\overline{g.\sigma}} \otimes 1.
\end{align*}
By applying Lemma \ref{lem:vacci}, since $g.\overline{\sigma}=\overline{g.\sigma}$, we have that Equation (\ref{eq:comma}) holds.

We consider now the following diagram:
$$
\begin{CD}
@>d_{i+1}>> \mathcal{R}_i @> d_i  >>  \mathcal{R}_{i-1}  @>d_{i-1}>> \\
@. @V\tilde{\phi}_{i}VV @V\tilde{\phi}_{i-1}VV\\
@>\tilde{\delta}_{n-i-2}>> C^{n-i-2} \otimes \mathbb{C}@ >\tilde{\delta}_{n-i-1} >> C^{n-i-1} \otimes \mathbb{C} @>\tilde{\delta}_{n-i}>>
\end{CD}
$$
And we define the $\tilde{\delta}_i$'s as an extension of the homomorphisms $\delta_i$:
$$
\tilde{\delta}_i(e^{*}_{\sigma}\otimes 1)=\underset{\sigma \cup k \in K^{*}}{\sum_{k \notin \sigma}} \ \textnormal{sgn}(k, \sigma \cup k) e^{*}_{\sigma \cup k} \otimes 1.
$$

From Equation (\ref{eq:trova}) it follows that:
$$
\tilde{\phi}_{i-1} \circ d_i=\tilde{\delta}_{n-i-1} \circ \tilde{\phi}_i.
$$
Thus, we have that:
$$
\alpha_{i+1} \simeq \rho^{*^{n-i-3}} \ \otimes \ \textnormal{sgn} .
$$
 
\end{proof}

\begin{lem}\label{teo:unoduetre}
Let $K$ be a simplicial complex with ground set $V$. Then:
$$
\widetilde{H}_{i}(K) \simeq \widetilde{H}_{i+1}((2^{V},K),\mathbb{C}).
$$
Furthermore if we consider the representations of the group $G$ on the reduced homology spaces of $K$ and $(2^{V},K)$,
$$
\begin{array}{cccc} 
\rho_i:& G & \longrightarrow & \textnormal{GL}(\widetilde{H}_{i}(K,\mathbb{C}))
\end{array} \quad \hbox{and} \quad \begin{array}{cccc} 
\alpha_{i+1}:& G & \longrightarrow & \textnormal{GL}(\widetilde{H}_{i+1}((2^{V},K),\mathbb{C})),
\end{array}
$$
we have that
$$
\rho_i \simeq \alpha_{i+1}.
$$
\end{lem}
\begin{proof}
 The isomorphism follows from Theorem 23.3 of \cite{munkres}: we have the long exact sequence of the pair $(2^{V},K)$:
$$
\cdots \xrightarrow \ \widetilde{H}_{i+1}(2^{V}) \xrightarrow \ \widetilde{H}_{i+1}(2^{V},K) \xrightarrow \ \widetilde{H}_{i}(K) \xrightarrow \ \widetilde{H}_{i}(2^{V}) \xrightarrow \ \cdots
$$
Since $2^{V}$ is the full simplex the spaces $\widetilde{H}_{i+1}(2^{V})$ and $\widetilde{H}_{i}(2^{V})$ are zero. Hence, the sequence becomes:
$$
\cdots \xrightarrow \ 0 \xrightarrow \ \widetilde{H}_{i+1}(2^{V},K) \xrightarrow \ \widetilde{H}_{i}(K) \xrightarrow \ 0 \xrightarrow \ \cdots
$$
It follows that the groups $\widetilde{H}_{i+1}(2^{V},K)$ and $\widetilde{H}_{i}(K)$ are isomorphic.

We now consider the following diagram:
$$
\begin{CD}
 \tilde{H}_{i+1}(2^{V},K) @> \partial^{*}  >>  \tilde{H}_{i}(K)   \\
 @V\alpha_{i+1,g}VV @V\rho_{i,g}VV\\
\tilde{H}_{i+1}(2^V,K) @ >\partial^{*} >> \tilde{H}_{i}(K) 
\end{CD}
$$
where $\partial^{*}$ is the \textit{homology boundary isomorphism} (see \cite{munkres}, Lemma 24.1):
$$
\begin{CD}
@. C_{i+1}(2^{V}) @> \pi_{\#}  >>  \mathcal{R}_{i+1}(2^V,K)   \\
@. @VV\partial^{V}_{i+1} V\\
C_{i}(K)@>i_{\#}>> C_{i}(2^V)
\end{CD}
$$
The isomorphism $\partial^{*}$ is defined by a certain \textit{zig-zag} process: pull back via $\pi_{\#}$, apply $\partial^{V}_{i+1}$, and pull back via $i_{\#}$. For each $g \in G$ we consider the action on the chain groups of the full simplex, of $K$ and of $(2^V,K)$:
$$
\begin{array}{cccc} 
\tilde g_{\#,i}^{V}:& C_{i}(2^V) & \longrightarrow & C_{i}(2^V)\\
& [j_{1},\dots,j_{i+1}] & \longmapsto & [g(j_{1}),\dots,g(j_{i+1})],
\end{array}
$$
$$
\begin{array}{cccc} 
\tilde g_{\#,i}:& C_{i}(K) & \longrightarrow & C_{i}(K)\\
& [j_{1},\dots,j_{i+1}] & \longmapsto & [g(j_{1}),\dots,g(j_{i+1})],
\end{array}
$$
$$
\begin{array}{cccc} 
\tilde g_{\#,i}^{V,K}:& \mathcal{R}_{i}(2^V,K) & \longrightarrow & \mathcal{R}_{i}(2^V,K).\\
& [j_{1},\dots,j_{i+1}] & \longmapsto & [g(j_{1}),\dots,g(j_{i+1})].
\end{array}
$$
We have that:
$$
\left.\tilde g_{\#,i}^{V}\right|_{C_{i}(K)}=\tilde{g}_{\#,i}, \qquad 
\left.\tilde g_{\#,i}^{V}\right|_{\mathcal{R}_i(2^V,K)}=\tilde{g}^{V,K}_{\#,i}.
$$

We also know that each boundary operator commutes with $\tilde{g}_{\#,i}$, $\tilde{g}^{V}_{\#,i}$ and $\tilde{g}^{V,K}_{\#,i+1}$ from Lemma 12.1 of \cite{munkres}. Let $b \in \tilde{H}_{i+1}(2^{V},K)$. There exists an $a \in \mathcal{R}_{i+1}(2^{V},K)$ such that $b=a+\textnormal{Im}(d_{i+2})$. Therefore:
\begin{align*}
    \rho_{i,g}(\partial^{*}(b))&=\rho_{i,g}(\partial^{V}_{i+1}(a)+\textnormal{Im}(d_{i+2}))=\tilde{g}_{\#,i}(\partial^{V}_{i+1}(a))+\textnormal{Im}(d_{i+2})\\
    &=\tilde{g}^V_{\#,i}(\partial^{V}_{i+1}(a))+\textnormal{Im}(d_{i+2})=
    \partial^{V}_{i+1}(\tilde{g}^V_{\#,i+1}(a))+\textnormal{Im}(d_{i+2}), \hbox{ and} \\
    \partial^{*}(\alpha_{i+1,g}(b))&=\partial^{*}(\tilde{g}^{V,K}_{\#,i+1}(a)+\textnormal{Im}(d_{i+2}))=\partial^{V}_{i+1}(\tilde{g}^{V,K}_{\#,i+1}(a))+\textnormal{Im}(d_{i+2})\\
    &=\partial^{V}_{i+1}(\tilde{g}^V_{\#,i+1}(a))+\textnormal{Im}(d_{i+2}).
\end{align*}
Thus, we have that 
$$
\partial^{*} \circ \alpha_{i+1,g}= \rho_{i,g} \circ \partial^{*} \quad \hbox{ for every } g \in G,
$$
and this implies that $\rho_i \simeq \alpha_{i+1}$.
\end{proof}
Combining the results of Lemma \ref{teo:unoduetre} and Lemma \ref{teo:trequattro} we obtain the proof of Theorem \ref{teo:fonda}.

\begin{remk}
From Alexander duality we know that for every simplicial complex $K$ on vertex set $V$ such that $V \notin K$, with $n=|V|$:

$$
\widetilde{H}_{i}(K)\simeq \widetilde{H}^{n-3-i}(K^{*}).
$$

In fact, working with complex coefficients the reduced cohomology group $\widetilde{H}^{j}(K)$ is the dual vector space of the reduced homology group $\widetilde{H}_{j}(K)$, so that 
$\widetilde{H}_{j}(K)\simeq \widetilde{H}^{j}(K) $. 
Combining the two results we obtain:
\begin{equation}\label{Deltino}
\widetilde{H}_{i}(K)\simeq \widetilde{H}_{n-3-i}(K^{*}).
\end{equation}
\end{remk}

\section{Equivariant cross-cut theory}\label{section:crosscut}
Let $L$ be a lattice with maximal and minimal elements $\widehat{1}$ and $\widehat{0}$ respectively.
    We recall the following definition from \cite{folk}:
\begin{defin}\label{def:crosscut}
If $L$ is a lattice with $\widehat{0}$ and $\widehat{1}$, a \textit{cross-cut} of $L$ is a set $C \subseteq L$ such that:
\begin{itemize}
\item[i)] $\widehat{0}$ , $\widehat{1} \notin C$.
\item[ii)] If $x,y \in C$ then $x \nless y $ and $y \nless x $. \qquad (\textit{$x$ and $y$ are incomparable})
\item[iii)] Any finite chain $x_{1}<x_{2}<\dots<x_{n}$ in $L$ can be extended to a chain which contains an element of $C$.
\end{itemize}
In particular, axiom iii) implies that every maximal chain contains an element of C.
\end{defin}

Let $L$ be a lattice with $\widehat{0}$ and $\widehat{1}$ and let $C$ be a \textit{cross-cut} of $L$.
\begin{defin}\label{def:span}
A finite subset $\{x_1, \dots, x_{n} \} \subseteq C$ \textit{\lq spans\rq} if and only if 
$$
x_{1}   \wedge x_{2} \wedge \dots \wedge x_{n}=\widehat{0} \qquad \hbox{and} \qquad x_{1}   \vee x_{2} \vee \dots  \vee x_{n}=
\widehat{1}
$$
\end{defin}
Here $x \wedge y$ denotes the largest element $\leq x$ and $\leq y$, and $x \vee y$ denotes the smallest element $\geq x$ and $\geq y$. 

Let $K(C)$ be the abstract simplicial complex whose vertices are the elements of $C$ and whose simplices are all finite subsets of $C$ which do not \textit{\lq span\rq}.
We denote $\widetilde{H}_{i}(C)=\widetilde{H}_{i}(K(C))$. 
Let $K(L)$ be the order complex of the lattice $L$ and denote $\widetilde{H}_{i}(L)=\widetilde{H}_{i}(K(L))$.
The following result was proved in \cite{folk}, Theorem 3.1:
\begin{teo}\label{folk1}
Let $L$ be a lattice and let $C$ be a cross-cut of $L$, then: $$\widetilde{H}_{i}(C) \simeq \widetilde{H}_{i}(L).$$
\end{teo}

In order to see that the previous isomorphism is also a $\mathbb{C}[G]$-module isomorphism we need the following result:
\begin{lem}\label{lem:barycentric}
Let $K$ be an abstract simplicial complex and let $K'$ be its first barycentric subdivision. Let also $G$ be a finite group of automorphisms of the face poset of $K$. Then we have the following isomorphism of $\mathbb{C}[G]$-modules:
$$
\widetilde{H}_{i}(K)\simeq_{G} \widetilde{H}_{i}(K').
$$
\end{lem}
\begin{proof}
First, we need to describe the action of $G$ on $K'$. Let $\mathcal{L}(K)$ be the face poset of $K$; it is clear that the order complex of $\mathcal{L}(K)$ is the barycentric subdivision of $K$. Thus, we have a straightforward $G$-action on the order complex of $\mathcal{L}(K)$ and its homology spaces. We have to show that the following two representations are isomorphic:
$$
\begin{array}{cccc} 
\tilde{\rho}_{i}:& G & \longrightarrow & \textnormal{GL}(H_i(K)) \\
& g & \longmapsto & \tilde{\rho}_{i,g} 
\end{array} \quad \hbox{and} \quad \begin{array}{cccc} 
\tilde{\rho}'_{i}:& G & \longrightarrow &\textnormal{GL}(H_i(K')). \\
& g & \longmapsto & \tilde{\rho}'_{i,g}
\end{array}
$$
Let $w*K$ be a cone. If $e_{\sigma}=[a_0,\dots,a_i]$ is an oriented simplex of $K$, let
$$
\Big[ w, e_{\sigma}\Big]=[w,a_0,\dots,a_i]
$$
denote an oriented simplex of $w*K$. This operation is well defined
and is called the \textit{bracket operation} (see \cite{munkres}, Section \S 8).

If $\sigma=\{a_0,\dots,a_i\}$ is a simplex, let $\hat{\sigma}$ denote the barycenter of $\sigma$. The complex $K'$ equals the collection of all simplices of the form
$$
[\hat{\sigma}_1,\dots,\hat{\sigma}_n] \quad \hbox{ where } \sigma_1\supset \cdots \supset \sigma_n.
$$
We know from \cite{munkres}, Section \S 17 that there is a unique augmentation-preserving chain map $\begin{array}{cccc} 
\textnormal{sd}:& C_i(K) & \longrightarrow & C_i(K') 
\end{array}$ called the \textit{barycentric subdivision operator} that induces an isomorphism of homology spaces. 
There is an inductive formula for the operator $\textnormal{sd}$. It is the following:
$$
\textnormal{sd}(v)=\hat{v}=v \qquad \hbox{for } v \in V, \hbox{ and}
$$
$$
\textnormal{sd}(e_{\sigma})=\Big[\hat{\sigma}, \textnormal{sd}(\partial_{i}(e_{\sigma})) \Big] \qquad \hbox{for } \sigma \in K \hbox{ with } \textnormal{dim}(\sigma)=i.
$$
Now we consider the following two representations:
$$
\begin{array}{cccc} 
\rho_{i}:& G & \longrightarrow & \textnormal{GL}(C_i(K)) \\
& g & \longmapsto & \rho_{i,g} 
\end{array} \quad \begin{array}{cccc} 
\rho'_{i}:& G & \longrightarrow &\textnormal{GL}(C_i(K')) \\
& g & \longmapsto & \rho'_{i,g}
\end{array}
$$
We want to show that the following diagram 
\begin{equation}\label{eq:diagrammino}
\begin{CD}
 C_i(K) @> \textnormal{sd}  >>  C_{i}(K')   \\
 @V\rho_{i,g}VV @V\rho'_{i,g}VV\\
C_i(K) @ >\textnormal{sd} >> C_{i}(K') 
\end{CD}
\end{equation}
commutes for every $g \in G$. We proceed by induction on $i$:
\begin{itemize}
    \item[-] Suppose $i=0$. It follows from the action of $G$ on the vertices of $K$ and $K'$ that $\rho_{0,g}(v)=\rho'_{0,g}(v)$ for every $v \in V$. Thus:
    $$
    \rho'_{0,g}(\textnormal{sd}(v))=\rho'_{0,g}(v)=\rho_{0,g}(v)=\textnormal{sd}(\rho_{0,g}(v)).
    $$
    \item[-] We now suppose the diagram commutes for $i=n$ and we prove it for $i=n+1$. Let $g.\sigma=\tau$, thus:
    \begin{align*}
        \rho'_{i+1,g}(\textnormal{sd}(e_{\sigma}))&=\rho'_{i+1,g}(\Big[\hat{\sigma}, \textnormal{sd}(\partial_{i+1}(e_{\sigma})) \Big])=\Big[\hat{\tau}, \rho'_{i,g}(\textnormal{sd}(\partial_{i+1}(e_{\sigma}))) \Big]\\
        &=\Big[\hat{\tau}, \textnormal{sd}(\rho_{i,g}(\partial_{i+1}(e_{\sigma}))) \Big]=\Big[\hat{\tau},\textnormal{sd}(\partial_{i+1}(\rho_{i+1,g}(e_{\sigma}))) \Big]\\
        &=\Big[\hat{\tau}, \textnormal{sd}(\partial_{i+1}(g.e_{\sigma})) \Big]=\textnormal{sd}(g.e_{\sigma})=\textnormal{sd}(\rho_{i+1,g}(e_{\sigma})).
    \end{align*}
\end{itemize}
Since the diagram (\ref{eq:diagrammino}) commutes and both $\textnormal{sd}$, $\rho_{i,g}$, $\rho'_{i,g}$ commute with the border operator $\partial$ we have that the following diagram commutes and consequently the lemma is proved:
$$
\begin{CD}
 \tilde{H}_i(K) @> \textnormal{sd}^{*}  >>  \tilde{H}_{i}(K')   \\
 @V\tilde{\rho}_{i,g}VV @V\tilde{\rho}'_{i,g}VV\\
\tilde{H}_i(K) @ >\textnormal{sd}^{*} >> \tilde{H}_{i}(K').
\end{CD}
$$
\end{proof}

\begin{defin}
Let $L$ be a lattice and $G$ a group of automorphism of $L$. A cross-cut $C$ of $L$ is \emph{$G$-stable} if $G.C=C$, i.e., if $C$ is the union of $G$-orbits.
\end{defin}

\begin{teo}\label{teo:crosscutfixed}
Let $L$ be a lattice and $G$ a group of automorphism of $L$. Let $C$ be a $G$-stable cross-cut of $L$. Then we have the following $\mathbb{C}[G]$-module isomorphism:
$$\widetilde{H}_{i}(L) \simeq_{G} \widetilde{H}_{i}(C).$$
\end{teo}
\begin{proof}
We briefly recall Folkman's argument. Let $K=K(L)$ be the order complex of $L$ and let $C=\{\alpha_1,\dots, \alpha_n\}$ be a cross-cut of $L$ fixed by $G$. For each $\alpha \in C$ let $L_{\alpha}$ be the subcomplex of $K$ consisting of all simplices $\{y_1,\dots,y_t\}$ such that the set $\{y_1,\dots,y_t,\alpha\}$ is totally ordered. By the third property of a cross-cut, the family $\{L_{\alpha}\}_{\alpha \in C}$ is a covering of $K$. In the proof of Theorem 3.1 (\cite{folk}) Folkman shows that $L_{\alpha_1}\bigcap \cdots \bigcap L_{\alpha_n}$ has the homology of a point or is empty and shows also that
\begin{equation}\label{eq:folkmanproof}
K(C)=\mathcal{N}(\{L_{\alpha}\}_{\alpha \in C})
\end{equation}
where $K(C)$ is the simplicial complex associated to the cross-cut $C$ and $\mathcal{N}=\mathcal{N}(\{L_{\alpha}\}_{\alpha \in C})$ is the nerve of the covering $\{L_{\alpha}\}_{\alpha \in C}$.
Thus, we can apply a nerve theorem. We follow the construction made by Bj\"orner in \cite{Bj4} Theorem 10.6. Let $P(K)$ and $P(\mathcal{N})$ be the face lattice associated to $K$ and $\mathcal{N}$, respectively. Bj\"orner defines the following order-reversing map of posets:
$$
\begin{array}{cccc} 
\tilde{f}:& P(K) & \longrightarrow & P(\mathcal{N}) \\
& \sigma & \longmapsto & \{\alpha \in C \ | \ \sigma \in L_{\alpha}\}.
\end{array}
$$
This map $\tilde{f}$ induces a simplicial map $f$ between the respective order complex of $P(K)$ and $P(\mathcal{N})$ which are the first barycentric subdivision of $K$ and $\mathcal{N}$:
$$
\begin{array}{cccc} 
f:& K' & \longrightarrow & \mathcal{N}' \\
& \{\sigma_0,\dots,\sigma_i\} & \longmapsto & \{\tilde{f}(\sigma_0),\dots,\tilde{f}(\sigma_i)\} 
\end{array}
$$
where $\sigma_0,\dots,\sigma_i$ are simplices of $K$, with $\sigma_0\supseteq \cdots \supseteq \sigma_i$ so that $\{ \sigma_0,\dots,\sigma_i\}$ is a simplex of $K'$.
Applying Theorem 10.6 of \cite{Bj4} we get, in particular, that $f$ induces a chain map $f_{\#}$ between $C_i(K')$ and $C_i(\mathcal{N}')$ in the following manner:
$$
f_{\#}([\sigma_0, \dots, \sigma_{i}])=\begin{cases} [\tilde{f}(\sigma_0), \dots, \tilde{f}(\sigma_{i})], & \mbox{if }\tilde{f}(v_0), \dots, \tilde{f}(v_i) \mbox{ are distinct} \\ 0, & \mbox{otherwise}
\end{cases}
$$
and moreover an isomorphism $f_{*}$ on homology spaces: 
$$
\widetilde{H}_{i}(K') \simeq \widetilde{H}_{i}(\mathcal{N}').
$$
We need to describe the action of $G$ on $K'$ and $\mathcal{N}'$: the $G$-action on $L$ induces an action on $K$ and therefore on $K'$ (in the sense of \ref{lem:barycentric}). Since $C$ is $G$-stable, every $g \in G$ acts on $C$ permuting its elements. Furthermore, since $g$ is an order automorphism of $L$, it acts on the covering $\{L_{\alpha}\}_{\alpha \in C}$ respecting the intersection relations. Therefore $G$ yields an action on the nerve $\mathcal{N}$ and therefore on $\mathcal{N}'$. We want to show that the following two representations are isomorphic:
$$
\begin{array}{cccc} 
\tilde{\rho}_{1}:& G & \longrightarrow & \textnormal{GL}(\tilde{H}_i(K')), \\
& g & \longmapsto & \tilde{\rho}_{1,g} 
\end{array} \quad \begin{array}{cccc} 
\tilde{\rho}_{2}:& G & \longrightarrow &\textnormal{GL}(\tilde{H}_i(\mathcal{N}')). \\
& g & \longmapsto & \tilde{\rho}_{2,g}
\end{array}
$$
Let 
$$
\begin{array}{cccc} 
\rho_{1}:& G & \longrightarrow & \textnormal{GL}(C_i(K')), \\
& g & \longmapsto & {\rho}_{1,g} 
\end{array} \quad \begin{array}{cccc} 
{\rho}_{2}:& G & \longrightarrow &\textnormal{GL}(C_i(\mathcal{N}')). \\
& g & \longmapsto & {\rho}_{2,g}
\end{array}
$$
be the representations on the chain spaces. Since $g\in G$ is an order automorphism of $L$ we have the following: if $\tilde{f}(\sigma)=\{\alpha_{j_0},\dots,\alpha_{j_t}\}=\beta$ then $\tilde{f}(g.\sigma)=\{g.\alpha_{j_0},\dots,g.\alpha_{j_t}\}=g.\beta$. We explicitly describe the maps induced by $g \in G$ on the chain spaces:
$$
\begin{array}{cccc} 
\rho_{1,g}:& C_i(K') & \longrightarrow & C_i(K') \\
& [\sigma_0, \dots, \sigma_{i}] & \longmapsto & [g.\sigma_0, \dots, g.\sigma_{i}],
\end{array}
$$
$$
\begin{array}{cccc} 
\rho_{2,g}:& C_i(\mathcal{N}') & \longrightarrow & C_i(\mathcal{N}'). \\
& [\beta_{j_0},\dots,\beta_{j_i}] & \longmapsto & [g.\beta_{j_0},\dots,g.\beta_{j_i}]
\end{array}
$$
where $\beta_{j}$ are simplices of $\mathcal{N}$ satisfying $\beta_{j_0}\subseteq \cdots \subseteq \beta_{j_i}$. We want to show that the following diagram commutes:
$$
\begin{CD}
 C_i(K') @> f_{\#}>>  C_{i}(\mathcal{N}')   \\
 @V{\rho}_{1,g}VV @V{\rho}_{2,g}VV\\
C_i(K') @ >f_{\#}>> C_{i}(\mathcal{N}'),
\end{CD}
$$

i.e., $\rho_{2,g}(f_{\#}([\sigma_0, \dots, \sigma_{i}]))=[g.\beta_{j_0},\dots,g.\beta_{j_i}]=
(f_{\#}(\rho_{1,g}([\sigma_0, \dots, \sigma_{i}])).$

Therefore the diagram commutes and since $f_{\#}$ is a chain map we have that $\tilde{\rho}_1 \simeq \tilde{\rho}_2$, i.e., $\widetilde{H}_{i}(K') \simeq_G \widetilde{H}_{i}(\mathcal{N}')$. Using the results of Lemma \ref{lem:barycentric} and Equation \ref{eq:folkmanproof} we have the following $\mathbb{C}[G]$-module isomorphism $$\widetilde{H}_{i}(K) \simeq_G \widetilde{H}_{i}((K(C))=\widetilde{H}_{i}(C).$$
\end{proof}

\begin{remk}
As one of the referees pointed out, in \cite{lakser} Lasker proved that $K(L)$ and $K(C)$ are homotopy equivalent. It could be shown that this homotopy equivalence is $G$-equivariant, which would imply another proof of Theorem \ref{teo:crosscutfixed}.
\end{remk}

\section{Applications to matroids}\label{section:matroid}

We now specialize the results of the previous two sections to matroids. For basic facts on matroids, the reader may refer to \cite{oxley}.
Let $M=(E,I)$ be a matroid with ground set $E$ and a collection of independent sets $I$, which forms an abstract simplicial complex. Let $M^{*}=(E,I^{*})$ be its dual. We recall that the rank of $A\subseteq E$ is the maximal cardinality of an element of $I$ contained in $A$. We say that
$A \subseteq E$ is \emph{non-spanning} in $M$ if $\textnormal{rk}(A)<\textnormal{rk}(E)$, i.e., $A$ does not contain any basis of $M$. 
Let
$$
NS(M)=\{A\subseteq E \ | \ A \hbox{ is non-spanning in } M\}.
$$

It is easy to see that $NS(M)$ is an abstract simplicial complex.

\begin{prop}\label{teononspan}
$A \subseteq E$ is non-spanning in $M^{*}$ if and only if $A^{c}$ is dependent in $M$.
\end{prop}
\begin{proof}
If $A \subseteq E$ is non-spanning in $M^{*}$ we have:
$$
\hbox{rk}^{*}(A)< \hbox{rk}^{*}(E).
$$
This is equivalent to:
$$
\hbox{rk}(A^{c}) + |A| - \hbox{rk}(E) < \hbox{rk}^{*}(E)
$$
and therefore to
$$
\hbox{rk}(A^{c})<- |A| +\hbox{rk}(E)+\hbox{rk}^{*}(E)=|E| - |A|=|A^{c}| \Longleftrightarrow 
A^{c} \notin I.
$$
For every $A \subseteq E$, we have $\textnormal{rk}(A)\leqslant |A|$, thus $A$ is independent if and only if $\textnormal{rk}(A)= |A|$.
\end{proof}

\begin{prop}
Let $\hbox{IN}(M)=I$ be the abstract simplicial complex associated with the independent sets of the matroid $M=(E,I)$ and let $I^{*}$ be its Alexander dual, then:
$$
I^{*}=NS(M^{*}).
$$
\end{prop}
\begin{proof}
Using the result shown in Proposition \ref{teononspan} we claim that:
\begin{align*}
I^{*}&=\{A \subseteq E  : A^{c} \notin I \}\\
&=\{A \subseteq E  : A^{c} \ \hbox{is dependent in} \ M=(E,I) \}\\
&=\{A \subseteq E  : A \ \hbox{is not spanning of} \ M^{*}\}=\hbox{NS}(M^{*}).
\end{align*}
\end{proof}
The previous result, together with Equation (\ref{Deltino}), implies the following:
$$
\widetilde{H}_{i}(\hbox{NS}(M))\simeq \widetilde{H}_{n-3-i}(\hbox{IN}(M^{*})).
$$
This is an isomorphism not only of vector spaces, but also of representations, up to a sign. Indeed, by applying Theorem \ref{teo:fonda}, we obtain:
\begin{teo}\label{teo:fond1}
Let $G$ be the automorphism group of a matroid. Then we have the following $\mathbb{C}[G]$-module isomorphism:
$$
\widetilde{H}_{i}(\hbox{NS}(M))\simeq_{G} \widetilde{H}_{n-3-i}(\hbox{IN}(M^{*})) \otimes \textnormal{sgn},
$$
where $n$ is the cardinality of the ground set of $M$. 
\end{teo}

Similarly, we can specialize the results from Section \ref{section:crosscut} to the case of matroids. Let $M=(E,I)$ be a simple matroid with $E=\{a_1,\dots a_m \}$. Let $\mathcal{L}(M)$ be the lattice of \emph{flats} of $M$ ordered by inclusion.
Since $M$ is simple, each singleton of $E$ is a flat. Thus $\{a_1\},\{a_2\}, \dots, \{a_m\} \in \mathcal{L}(M)$ and each corresponds to an atom of the poset $(\mathcal{L}(M), \subseteq)$. We now consider a set $C$ defined as
$$
C=\{\{a_1\},\{a_2\}, \dots, \{a_m\} \} \subseteq \mathcal{L}(M).
$$
Since $C$ satisfies the three axioms of Definition \ref{def:crosscut}, the set $C$ is a \textit{cross-cut} of $\mathcal{L}$.
We want to prove that:
$$
K(C)=NS(M).
$$
In the following proposition we perform a slight abuse of notation by identifying:
$$
C=\{\{a_1\},\{a_2\}, \dots, \{a_m\} \}=\{a_1,a_2, \dots, a_m\}. 
$$

\begin{prop}\label{teo:span}
$A \subseteq C$ does not \textit{\lq span\rq} (in the sense of Definition \ref{def:span}) if and only if $A$ is a non-spanning set in $M=(E,I)$.
\end{prop}
\begin{proof}
\leavevmode
\begin{itemize}
\item[$\Longrightarrow)$] In $\mathcal{L}(M)$ we have:
$$
\widehat{0}=\emptyset \qquad \hbox{and} \qquad \widehat{1}=E.
$$
Let $A=\{a_{i_{1}},a_{i_{2}},\dots,a_{i_{n}}\}$ be a subset of $C$. If $A \subseteq C$ does not \textit{\lq span\rq}:
\begin{equation}\label{villa}
a_{i_{1}}\vee a_{i_{2}} \vee \dots \vee a_{i_{n}}=D\neq \widehat{1} 
\end{equation}
$D \in \mathcal{L}(M)$ and $D \neq \widehat{1}$ implies that $D$ is a non-spanning subset of $E$ because the only spanning subset in $\mathcal{L}(M)$ is $E=\widehat{1}$.

It follows from (\ref{villa}) that $A \subseteq D$; since $D$ is a non-spanning subset of $E$ therefore $A$ is a non-spanning subset of $E$.
\item[$\Longleftarrow)$] In $NS(M)$ the bases are the maximal non-spanning subsets of E, (i.e., the subsets of $E$, such that if we add an element they become spanning set) so they are flats, in particular they correspond to the co-atoms of $(\mathcal{L}(M),\subseteq)$.

Let $A=\{a_{i_{1}},a_{i_{2}},\dots,a_{i_{n}}\}$ be a non-spanning subset of $E$, there exist a basis $\mathcal{B}$ of $NS(M)$ such that:
$$
\hbox{if }A \subseteq \mathcal{B} \hbox{ and } \mathcal{B} \hbox{ is a flat, then }  \mathcal{B} \in \mathcal{L}(M).
$$
This implies:
$$
a_{i_{1}}\vee a_{i_{2}} \vee \dots \vee a_{i_{n}} \subseteq \mathcal{B} \neq \widehat{1}
$$
therefore $A$ does not \textit{\lq span\rq}.
\end{itemize}
\end{proof}
Using the result of Proposition \ref{teo:span}, we obtain:
$$
K(C)=NS(M).
$$
Since $C$ is a cross-cut fixed by $G$, we can apply Theorem \ref{teo:crosscutfixed} to $\mathcal{L}(M)$ and $C$ itself:
\begin{teo}\label{teo:fond2}
Let $G$ be the group of automorphism of the simple matroid $M$. Then we have the following $\mathbb{C}[G]$-module isomorphism:
$$
\widetilde{H}_{i}(\mathcal{L}(M))\simeq_{G} \widetilde{H}_{i}(C)=\widetilde{H}_{i}(NS(M))
$$
where $C$ is the cross-cut of $\mathcal{L}(M)$ composed of its atoms.
\end{teo}

By combining Theorem \ref{teo:fond1} and Theorem \ref{teo:fond2}, we get the following theorem:
\begin{teo}\label{teo:combined}
Let $G$ be the group of automorphism of the simple matroid $M$. Then we have the following $\mathbb{C}[G]$-module isomorphism:
$$
\widetilde{H}_{n-3-i} (IN(M^{*})) \simeq_{G} \widetilde{H}_{i} (\mathcal{L}(M)) \otimes \textnormal{sgn}
$$
where $n$ is the cardinality of the ground set of $M$.
\end{teo}

\section{Top cohomology of hyperplane arrangements}\label{section:arrangement}
Let $\mathcal{A}$ be a central arrangement of hyperplanes in $\mathbb C ^r$ and let $\mathcal{L}(\mathcal{A})$ be its intersection lattice. Let $M_{\mathcal{A}}$ be the matroid associated with $\mathcal{A}$; then the lattice of flats $\mathcal{L}(M_{\mathcal{A}})$ of $M_{\mathcal{A}}$ is isomorphic to $\mathcal{L}(\mathcal{A})$. We can assume that the arrangement is essential: then the rank of the matroid is $r$. We define the complement of the arrangement:
$$
\mathcal{C}(\mathcal{A})=\mathbb{C}^{r} \smallsetminus 
\bigcup_{H \in \mathcal{A}} H .
$$
Let $G$ be a subgroup of $GL(\mathbb{C}^r)$ that permutes the elements of $\mathcal{A}$; it is easy to see that $G$ is also a group of automorphism of the matroid $M_{\mathcal{A}}$. Let $\mathfrak{A}$ be the Orlik-Solomon algebra associated to $\mathcal{L}(\mathcal{A})$, and let $\mathfrak{B}$ be the \textit{algebra defined by shuffle} defined respectively in Section 2 and Section 3 of \cite{orliksolomon}. These algebras are $\mathbb{Z}-$graded: we denote by $\mathfrak{A}_r$ and $\mathfrak{B}_r$ the direct summands corresponding to the top degree $r$. In Theorem 3.7 of the same paper, Orlik and Solomon provide a $G$-isomorphism:
$$
\begin{array}{cccc} 
\theta:& \mathfrak{A} & \longrightarrow & \mathfrak{B}.
\end{array}
$$
Furthermore, we state Theorem 4.3 of \cite{orliksolomon}:
\begin{teo}\label{teo:orliklattice}
Let $L$ be a finite geometric lattice of rank $r>1$. Then $\mathfrak{B}_{r}$ and $H_{r-2}(L)$ are isomorphic $\mathbb{C}[G]$-modules.
\end{teo}
Combining the previous results we get the following $\mathbb{C}[G]$-module isomorphism:
\begin{equation}\label{eq:hyperplanelattice}
H^{r}(\mathcal{C}(\mathcal{A})) \simeq_{G} \mathfrak{A}_{r} \simeq_G \mathfrak{B}_{r} \simeq_G      H_{r-2}(\mathcal{L}(\mathcal{A})).
\end{equation}

Applying Theorem \ref{teo:combined} we obtain the following:
\begin{teo}\label{teo:iperpiani}
Let $\mathcal{A}$ be a central essential hyperplane arrangement of dimension $r$ and let $M_{\mathcal{A}}$ be the associated matroid with ground set of cardinality $n$. Then we have the following  $\mathbb{C}[G]$-module isomorphism:
$$
H^{r}(\mathcal{C}(\mathcal{A}))  \simeq_{G} {H}_{n-r-1} (IN(M^{*}_{\mathcal{A}})) \otimes \textnormal{sgn}.
$$
\end{teo}
In \cite{lehrersolomon} Lehrer and Solomon conjecture that if $W$ is a Coxeter group and $\mathcal{A}_W$ is the hyperplane arrangement associated to $W$ then there is a $\mathbb{C}[G]$-module isomorphism
$$
H^{p}(\mathcal{C}(\mathcal{A}_W)) \simeq_{W} \bigoplus_c \textnormal{Ind}^{W}_{Z(c)}(\xi_c) \qquad p=0,\dots,\textnormal{rank}(W)
$$
where $c$ runs over a set of representatives for the conjugacy classes of $W$ such that the dimension of the image of $c$ (viewed as an element of $GL(V)$) is equal to $p$ and $\xi_c$ is a suitable character of the centralizer $Z(c)$ of $c$ in $W$. They proved the conjecture for group of rank 2 and for $W=\mathfrak{S}_r$. In the case of the symmetric group $\mathfrak{S}_r$ the arrangement $\mathcal{A}_{\mathfrak{S}_r}$ is the braid arrangement and the intersection lattice $\mathcal{L}(\mathcal{A})$ is the \textit{partition lattice} $\Pi_r$, that is, the family of all partitions of the set $\{1,\dots,r\}$ partially ordered by refinement. 
Stanley studied the representations on the homology of the partition lattice in \cite{stan1}. By Equation (\ref{eq:hyperplanelattice}), his result agrees with the conjecture of Lehrer and Solomon. 

We remark that Theorem \ref{teo:iperpiani} allows us to rewrite Lehrer and Solomon' conjecture in the top cohomology case in the language of matroids:
$$
{H}_{n-r-1} (IN(M^{*}_{\mathcal{A}_W})) \simeq_{W} \bigoplus_c \textnormal{Ind}^{W}_{Z(c)}(\xi_c)) \otimes \textnormal{sgn}.
$$

\section{Coned graphs and complete
bipartite graphs}\label{section:kook}
In \cite{kook}, Woong Kook studied the homology of the independence complex $IN(M( \widehat{\Gamma} ))$ of the matroid associated to a \emph{coned graph} $\widehat{\Gamma}$, i.e. the graph obtained by adding a new vertex $p$ to a graph $\Gamma$ and joining each vertex of $\Gamma$ to $p$ by a simple edge. We recall the following definition from (\cite{kook}, Section 2):
\begin{defin}
An \textit{edge-rooted forest} $(F,\textbf{e})$ in $\Gamma$ is a spanning forest $F$ that contains at least one edge for each connected component of $\Gamma$, together with the datum $\textbf{e}$ of one edge for each component (called \emph{edge root}). 
\end{defin}
The rank of the only non zero homology group of $IN(M( \widehat{\Gamma}))$ is shown to be equal to the cardinality of the set of \textit{edge-rooted forests} $\mathcal{F}_{e}(\Gamma)$ in $\Gamma$.

In Section 3, Kook constructs a basis $\{ z_{F,\textbf{e}} \ : \ (F, \textbf{e}) \in \mathcal{F}_e(\Gamma)\}$ for $\widetilde{H}_{n-1}(IN(M(\widehat{\Gamma})))$ (where $n$ is the number of vertices of $\Gamma$).
This basis is indexed by the elements $(F,\textbf{e}) \in \mathcal{F}_e(\Gamma)$.
In the same Section, Kook describes the action of the automorphism group $G =\hbox{Aut}(\Gamma)$ on $\widetilde{H}_{n-1}(IN(M(\widehat{\Gamma})))$ for a finite simple graph $\Gamma$, showing that this action is isomorphic to the permutation action on $\mathcal{F}_e(\Gamma)$ tensored with the sign representation:
$$
g(z_{F,\textbf{e}})=\hbox{sgn}(g) \ z_{g(F,\textbf{e})}
$$
where $g(F,\textbf{e})=(g(F),g(\textbf{e}))$ (see Theorem 6, \cite{kook}). Extending by linearity those two $G$-actions, we obtain two representations of $G$, respectively on $\widetilde{H}_{n-1}(IN(M(\widehat{\Gamma})),\mathbb{C})$ and on the vector space $\boldsymbol{{\mathcal{F}_e(\Gamma)}}$ of formal $\mathbb{C}$-linear combinations of elements of $\mathcal{F}_e(\Gamma)$, which are isomorphic up to a sign:
$$
\widetilde{H}_{n-1}(IN(M(\widehat{\Gamma}))) \simeq_{G } \boldsymbol{{\mathcal{F}_e(\Gamma)}} \otimes \hbox{sgn} .
$$

Applying Theorem \ref{teo:combined} we obtain the following $\mathbb{C}[G ]$-module isomorphism:
\begin{equation}\label{gra1}
\widetilde{H}_{l-n-2} (\mathcal{L}(M^*(\widehat{\Gamma}))) \simeq_{G } \boldsymbol{{\mathcal{F}_e(\Gamma)}}
\end{equation}
where $l$ is the number of edges of $\widehat{\Gamma}$.

Furthermore, in \cite{kook2} Woong Kook and Kang-Ju Lee studied the homology of the independence complex $IN(M( K_{m+1,n+1}))$ of the matroid associated to the complete bipartite graph $K_{m+1,n+1}$.
We need the following definition:
\begin{defin}
A \textit{B-edge-rooted forest} $(F,\textbf{b},\textbf{e})$ in a complete bipartite graph $K_{m,n} (m,n \geq 1)$ is a spanning forest $F$ in $K_{m,n}$ composed of two kinds of connected components such that
\begin{itemize}
    \item[-] exactly one component is \textit{bi-rooted}, i.e., has one vertex-root in each bipartite set;
    \item[-] each of the remaining components is edge-rooted, i.e., has one edge marked as edge-root. (See Definition 3.3 in \cite{kook2}).
\end{itemize}

\end{defin}
The rank of the only non zero homology group of $IN(M( K_{m+1,n+1}))$ is shown to be equal to the cardinality of the set of the \textit{B-edge-rooted forests} $\mathcal{F}_{e}^{B}(K_{m,n})$ in $K_{m,n}$. 

In Section 5, the authors construct a basis $\{ z_{F,\textbf{b},\textbf{e}} \ : \ (F,\textbf{b},\textbf{e}) \in \mathcal{F}_e^{B}(K_{m,n})\}$ for $\widetilde{H}_{m+n}(IN(M( K_{m+1,n+1})))$. This basis is indexed by the elements $(F,\textbf{b},\textbf{e}) \in \mathcal{F}_e^{B}(K_{m,n})$.
In the same section, they proved the following theorem:
\begin{teo}
The action of $\mathfrak{S}_m \times \mathfrak{S}_n$ as a subgroup of $\mathfrak{S}_{m+1} \times \mathfrak{S}_{n+1}$ on\\ $\widetilde{H}_{m+n}(IN(M( K_{m+1,n+1})))$ is isomorphic to the action on $\mathcal{F}_e^{B}(K_{m,n})$ tensored with the sign representation:
    \[\sigma(z_{F,\textbf{b},\textbf{e}})=\textnormal{sgn}(\sigma) z_{\sigma(F,\textbf{b},\textbf{e})}.\]
\end{teo}
Now we consider the representations of the group $\mathfrak{S}_m \times \mathfrak{S}_n$ that extend by linearity the two $\mathfrak{S}_m \times \mathfrak{S}_n$- actions, respectively on the vector space $\widetilde{H}_{m+n}(IN(M( K_{m+1,n+1})))$ and on the vector space $\boldsymbol{{\mathcal{F}_e^{B}(K_{m,n})}}$ of formal $\mathbb{C}$-linear combinations of elements of $\mathcal{F}_e^{B}(K_{m,n})$. Clearly we have the following $\mathbb{C}[\mathfrak{S}_m \times \mathfrak{S}_n]$-module isomorphism:
$$
\widetilde{H}_{m+n}(IN(M( K_{m+1,n+1}))) \simeq_{\mathfrak{S}_m \times \mathfrak{S}_n} \boldsymbol{{\mathcal{F}_e^{B}(K_{m,n})}} \otimes \hbox{sgn} .
$$
Applying Theorem \ref{teo:combined} we obtain the following $\mathbb{C}[\mathfrak{S}_m \times \mathfrak{S}_n]$-module isomorphism:
\begin{equation}\label{gra2}
\widetilde{H}_{l-m-n-3} (\mathcal{L}(M^*(K_{m+1,n+1}))) \simeq_{\mathfrak{S}_m \times \mathfrak{S}_n} \boldsymbol{{\mathcal{F}_e^{B}(K_{m,n})}}
\end{equation}
where $l$ is the number of edges of $K_{m+1,n+1}$.

\section{The dual matroid of the complete graph}\label{section:completegraph}
We now consider the matroid $M(K_m)$ of the complete graph $K_m$, which has rank $r=m-1$ and ground set of cardinality $n=\binom{m}{2}$. This matroid is isomorphic to the matroid $M(\Phi^{+}_{A_{m-1}})$ associated with the positive roots of the root system of type $A_{m-1}$. In fact, this is the case of interest in \cite{DHM}.

We recall that the lattice of flats of this matroid is isomorphic to the {partition lattice} $\Pi_m$. 
In this case, Theorem \ref{teo:combined} specializes to the following:

\begin{teo}\label{teo:supremofinale}
$\widetilde{H}_{n-3-i} (IN(M^{*}(K_m)))$ and $\widetilde{H}_{i}(\Pi_m) \otimes \textnormal{sgn}$ are isomorphic as $\mathfrak{S}_m$-modules for every $i\ge0$.
\end{teo}
\begin{remk}\label{miglio}
In \cite{DHM}, de Cataldo, Heinloth and Migliorini apply this result to the computation of the supports of the perverse cohomology sheaves of the Hitchin fibration for $GL_m$ over the locus of reduced spectral curves.
\end{remk}

Rephrased in terms of root system of type $A_{m-1}$, the theorem above yields the following $\mathbb{C}[G]$-module isomorphism:
\begin{equation}\label{eq:rootsystem}
\widetilde{H}_{n-3-i} (IN(M^{*}(\Phi^{+}_{A_{m-1}}))) \simeq_{\mathfrak{S}_m} \widetilde{H}_{i}(\Pi_m) \otimes \textnormal{sgn}
\end{equation}
where
$$
n=|E(M^{*}(\Phi^{+}_{A_{m-1}}))|=|\Phi^{+}(A_{m-1})|=\binom{m}{2}=\frac{m(m-1)}{2}.
$$

\begin{remk}
We can make a dimensional calculation to better understand the dimensional shift.
The matroid $M(\Phi^{+}_{A_{m-1}},I)$ has rank equal to $m-1$, i.e. each basis has $m-1$ elements. Therefore, the matroid 
$M^{*}(\Phi^{+}_{A_{m-1}},I)$ has rank equal to:
$$
n-(m-1)=\frac{m(m-1)}{2}-(m-1)=\frac{(m-1)(m-2)}{2}.
$$
Thus, the dimension of the top homology of $IN(M^{*}(\Phi^{+}_{A_{m-1}})$ is one less than the number of the elements of a basis of $M^{*}(\Phi^{+}_{A_{m-1}},I)$:
$$
\frac{(m-1)(m-2)}{2}-1.
$$
By Equation \ref{eq:rootsystem} we have the following isomorphism of $\mathbb{C}$-vector spaces:
$$
\widetilde{H}_{n-3-i} (IN(M^{*}(\Phi^{+}_{A_{m-1}}))) \simeq \widetilde{H}_{i}(\Pi_m).
$$
We impose
$$
n-3-i=\frac{(m-1)(m-2)}{2}-1
$$
then we have $i=m-3$ from $n=m(m-1)/2$.
Indeed, $H_{m-3}(\Pi_m)$ is the only nonzero homology group of $\Pi_m$. 
\end{remk}
By Theorem \ref{teo:supremofinale} these two representations
$$ 
\begin{array}{cccc} 
\rho_{n-m}: & \mathfrak{S}_m & \longrightarrow &  \hbox{GL}(\widetilde{H}_{n-m} (IN(M^{*}(\Phi^{+}_{A_{m-1}})))) \\
\end{array}
$$
and
$$
\begin{array}{cccc} 
\gamma_{m-3}: & \mathfrak{S}_m & \longrightarrow &  \hbox{GL}(\widetilde{H}_{m-3}(\Pi_m)\otimes \textnormal{sgn}) \\
\end{array}
$$
are isomorphic. From a result due to Stanley (\cite{stan1}, Theorem 7.3) we know that the representations on the top homology of the partition lattice
$$
\begin{array}{cccc} 
\tilde{\gamma}_{m-3}: & \mathfrak{S}_m & \longrightarrow &  \hbox{GL}(\widetilde{H}_{m-3}(\Pi_m)) \\
\end{array}
$$
are the following 
$$
\tilde{\gamma}_{m-3} \simeq \textnormal{sgn} \otimes \textnormal{ind}_{C_m}^{\mathfrak{S}_m}(e^{2\pi i / m}).
$$

Thus, we get 
$$
\rho_{n-m} \simeq \textnormal{ind}_{C_m}^{\mathfrak{S}_m}(e^{2\pi i / m}) 
$$
or as $\mathbb{C}[\mathfrak{S}_m]$-modules:
$$
\widetilde{H}_{n-m} (IN(M^{*}(\Phi^{+}_{A_{m-1}}))) \simeq_{\mathfrak{S}_m} \textnormal{ind}_{C_m}^{\mathfrak{S}_m}(e^{2\pi i / m}).
$$


\newpage
\begin{thebibliography}{9}
        \bibitem{Bac1}\label{Bac1} {\sc K. Baclawsky and A. Bj$\ddot{\hbox{o}}$rner}, \textit{Fixed points in partially ordered set}, Advances in Math. \textbf{31} (1979), 263-287.
        \bibitem{Bj1}\label{Bj1}{\sc A. Bj$\ddot{\hbox{o}}$rner}, \textit{Homology and shellability of geometric lattices}, Matroid Applications, ed. N. White, Cambridge University Press, 1992, pp. 226-283.
        \bibitem{Bj3}\label{Bj3}{\sc A. Bj$\ddot{\hbox{o}}$rner}, \textit{On the homology of geometric lattices}, Algebra Universalis \textbf{14} (1982),107-128.
        \bibitem{Bj4}\label{Bj4}{\sc A. Bj$\ddot{\hbox{o}}$rner}, \textit{Topological Methods}, Handbook of Combinatorics, ed. R. Graham, M. Gr\"otschel, L. Lov\'{a}sz, North Holland, Amsterdam, 1995, 1819-1872.
        \bibitem{Bj2}\label{Bj2}{\sc A. Bj$\ddot{\hbox{o}}$rner and M. Tancer}, \textit{Note: Combinatorial Alexander Duality - A Short and Elementary Proof}, Discrete \& Computational Geometry \textbf{42} (2009), 586-593.
        \bibitem{DHM} { \sc M. de Cataldo, J. Heinloth and L. Migliorini}, \textit{A support theorem for the Hitchin fibration: the case of $GL_n$ and $K_C$}, arXiv:1906.09582.
        \bibitem{folk}\label{folk} {\sc J. Folkman}, \textit{The Homology Groups of A Lattice}, J. Math. Mech. \textbf{15} (1966), 631-636.
        \bibitem{han}\label{han} {\sc P. Hanlon}, \textit{The fixed-point partition lattices}, Pacific J. Math. \textbf{96} (1981), 319-341.
        \bibitem{hatcher}\label{hatcher} {\sc A. Hatcher}, \textit{Algebraic Topology}, Cambridge University Press, New York (2010).
        \bibitem{hump}\label{hump} {\sc J.E. Humphreys}, \textit{Introduction to Lie Algebras and Representation Theory}, Springer, New York (1975).
        \bibitem{kook}{\sc W. Kook}, \textit{The homology of the cycle matroid of a coned graph}, European J. Combin. \textbf{28} (2007), 734-741.
        \bibitem{kook2}{\sc W. Kook and K.J. Lee}, \textit{Mobius coinvariants and bipartite edge-rooted forests}, European J. Combin. \textbf{71} (2018), 180-193.
        \bibitem{lakser}{\sc H. Lakser}, \textit{The homology of a lattice}, Discrete Math. \textbf{1} (1971), 187-192.
        \bibitem{lehrersolomon}{\sc G.I. Lehrer and L. Solomon}, \textit{On the action of the symmetric group on the cohomology of the complement of its reflecting hyperplanes}, J. Algebra \textbf{104} (1986), 410-424.
        \bibitem{munkres}\label{munkres}{ \sc J.R. Munkres}, \textit{Elements of Algebraic Topology}, The Benjamin/Cummings Publishing Company, Inc., Menlo Park (1984).
        \bibitem{orliksolomon}{\sc P. Orlik and L. Solomon}, \textit{Combinatorics and topology of complements of hyperplanes}, Invent. Math. \textbf{56} (1980), 167-189. 
        \bibitem{oxley}\label{oxley} {\sc J.G. Oxley}, \textit{Matroid Theory}, Oxford Univ. Press, (1992).
        \bibitem{rota1}\label{rota1} {\sc G.C. Rota}, \textit{On the foundation of combinatorial theory. I. Theory of Mobius functions}, Z. Wahrscheinlichkeitstheorie und Verw. Gebiete \textbf{2} (1964), 340-368.
        \bibitem{serre}\label{serre}{\sc J.P. Serre}, \textit{Linear Representations of Finite Groups}, Springer-Verlag, New York Berlin Heidelberg (1977).
        \bibitem{stan1}\label{stan1}{\sc R.P. Stanley}, \textit{Some aspects of groups acting on finite posets}, J. Combin. Theory (A) \textbf{32} (1982), 132-161.
         
\end{thebibliography}
\end{document}